\documentclass[12pt]{article}
\usepackage{amssymb,amsthm,amsmath}
\usepackage{setspace}
\usepackage{bbm}
\def\bm#1{\mathbbm{#1}}
\def\r#1{{\rm #1}}
\def\fn#1{\mathop{{\rm #1}\vphantom{\dim}}}


\newtheorem{theorem}{\indent Theorem}[section]

\newtheorem{defn}[theorem]{\indent Definition}
\newtheorem{lemma}[theorem]{\indent Lemma}
\newtheorem{prop}[theorem]{\indent Proposition}
\newtheorem{coro}[theorem]{\indent Corollary}
\newtheorem{remark}{\indent Remark}[section]

\renewcommand{\proofname}{\indent Proof.}

\title{Essential dimensions of finite groups}
\author{{ Ming-chang Kang}\\[2mm]
Department of Mathematics\\
National Taiwan University\\
Taipei, Taiwan\\
E-mail: kang@math.ntu.edu.tw}
\date{}

\begin{document}

\maketitle
\footnote{\hspace*{-7.5mm} 2000 Mathematics Subject Classification: 12E05, 12F10, 12F20, 13F30.\\
Keywords and Phrases: Essential dimension, Galois theory, group actions, valuation rings.}

\noindent Abstract. Let $K$ be an arbitrary field and $G$ be a
finite group. We will study the essential dimension of $G$ over
$K$, which is denoted by $\fn{ed}_K(G)$. A generalization of the
central extension theorem of Buhler and Reichstein (Compositio \
Math.\ {\bf 106} (1997) 159--179, Theorem 5.3) is obtained. As a
corollary, it can be proved that $\fn{ed}_K(S_n)\ge \lfloor
\frac{n}{2}\rfloor$ and $\fn{ed}_K(A_n)\ge 2\lfloor
\frac{n}{4}\rfloor$ for any field $K$ with $\fn{char} K\ne 2$,
while $\fn{ed}_K(S_n)\ge \lfloor\frac{n+1}{3}\rfloor$,
$\fn{ed}_K(A_n)\ge\lfloor\frac{n}{3}\rfloor$ for any field $K$
with $\fn{char} K=2$.

\newpage
\section{Introduction}

Let $G$ be a finite group and $K$ be any field (it is unnecessary
to assume that $K$ is algebraically closed or $\fn{char} K=0$).
The notion of the essential dimension of $G$ over $K$, denoted by
$\fn{ed}_K(G)$, was introduced by Buhler and Reichstein in
\cite{BR} when $\fn{char}K=0$. This concept was generalized for an
arbitrary field $K$ by Jensen, Ledet and Yui \cite{JLY}. Similar
notions were studied by Reichstein when $G$ was an algebraic group
\cite{Re}. Subsequently a functorial approach was proposed by
Merkurjev \cite{Me}, which also provided the definition of
$\fn{ed}_K(G)$ for an arbitrary field $K$. Merkurjev's method was
investigated in details by Berhuy and Favi \cite{BF}. As an
application of the functorial approach, Berhuy and Favi were able
to determine the essential dimensions of cyclic and dihedral
groups over various fields \cite[Section 7]{BF}.

The purpose of this article is to generalize most results of
Buhler and Reichstein in \cite{BR} without assuming
$\fn{char}K=0$. Instead of the functorial method, we follow the
method of \cite{BR} and \cite{JLY}. For example, one of the main
results is a generalization of the central extension theorem of
Buhler and Reichstein \cite[Theorem 5.3]{BR} (see Theorems
\ref{t4.5} and \ref{t4.6}). As applications of it, we are able to
prove the following theorems.

\begin{theorem} \label{t1.1}
Let $p$ be a prime number and $K$ be any field with $\fn{char}K\ne
p$ such that $\zeta_p\in K$. Then
$\fn{ed}_K((\bm{Z}/p\bm{Z})^r)=r$.
\end{theorem}

\begin{theorem} \label{t1.2}
Let $K$ be an arbitrary field and $S_n$ be the symmetric group.

\r{(1)} If $n\ge 5$, then $\fn{ed}_K(S_n)\le n-3$.

\r{(2)} $\fn{ed}_K(S_2)=\fn{ed}_K(S_3)=1$ and $\fn{ed}_K(S_4)=\fn{ed}_K(S_5)=2$.

\r{(3)} If $\fn{char}K\ne 2$, then $\fn{ed}_K(S_6)=3$.
\end{theorem}

Note that Theorem \ref{t1.1} and Part (2) and Part (3) of Theorem
\ref{t1.2} were proved by Berhuy and Favi \cite[Corollary
4.16]{BF} and Kraft \cite{Kr} respectively by different methods.
Part (1) of Theorem \ref{t1.2} was proved in \cite{BR} when
$\fn{char}K=0$, and was proved in \cite[p.318]{BF} when
$\fn{char}K\ne 2$.

We will organize this paper as follows. For the convenience of the
reader, we will recall the basic notions of essential dimensions
in Section 2, following the approach of \cite{BR} and \cite{JLY}.
In Section 3 we prove that the field of cross-ratios with respect
to the $S_n$ action over an arbitrary field $K$ is purely
transcendental of degree $n-3$ over $K$ (see Lemma \ref{l3.1}).
This result was proved in \cite{BR} when $\fn{char}K=0$. Part (1)
of the above Theorem \ref{t1.2} will be proved (see Proposition
\ref{p3.2}). The main result of this paper, i.e.\ the central
extension theorem, will be proved in Section 4. See Theorem
\ref{t4.5} and Theorem \ref{t4.6}. The proof of Theorem \ref{t1.1}
will be given in Section 4 (see Theorem \ref{t4.7}). In Section 5,
the essential dimensions of symmetric and alternating groups will
be studied (see Theorem \ref{t5.4}, Theorem \ref{t5.6} and the
proof of Theorem \ref{t1.2}). If $K$ is an arbitrary field, we
will give a necessary and sufficient condition for
$\fn{ed}_K(D_n)=1$ in Theorem \ref{t5.8} where $D_n$ is the
dihedral group of order $2n$. We will prove a characterization of
finite groups $G$ with $\fn{ed}_K(G)=1$ in a forthcoming article
\cite{CHKZ}.

Standing notation and terminology. For emphasis, $K$ is an
arbitrary field. All the fields in this article are assumed to
contain the ground field $K$. If $E$ is a field extension of $K$,
$\fn{trdeg}_K E$ denotes the transcendence degree of $E$ over $K$.
$\bm{F}_q$ denotes the finite field with $q$ elements. The order
of an element $\sigma$ in a group $G$ is denoted by
$\fn{ord}(\sigma)$. All the groups $G$ in this article are finite
groups. When we say that $G\to GL(V)$ is a representation of a
finite group $G$, it is understood that $V$ is a
finite-dimensional vector space over $K$. We will adopt the
following notations,

\begin{quotation}
$S_n$, the symmetric group,

$A_n$, the alternating group,

$\bm{Z}/n\bm{Z}$, the cyclic group of order $n$,

$PGL_2(K)$, the group isomorphic to $GL_2(K)/K^{\times}$,

$D_n$, the dihedral group of order $2n$.
\end{quotation}

We will take the convention that $\fn{char}K\nmid  n$ means either
$\fn{char}K=0$ or $\fn{char}K=p>0$ with $p\nmid n$. Finally
$\zeta_n$ denotes a primitive $n$-th root of unity; whenever we
write $\zeta_n \in K$, it is assumed tacitly that $\fn{char}K\nmid
n$.

\section{Definitions and preliminaries}

Let $K$ be an arbitrary field.

\begin{defn}\label{d2.1} \rm
Let $L/L_0$ be a finite separable extension of fields containing the field $K$.
$L/L_0$ is defined over $E_0$ if there is a field $E$ such that $L\supset E\supset E_0 \supset K$,
$[E:E_0]=[L:L_0]$ and $L=EL_0$.
The essential dimension of $L/L_0$, denoted by $\fn{ed}_K(L/L_0)$,
is the minimum of $\fn{trdeg}_KE_0$,
where $E_0$ runs through all intermediate fields over which $L/L_0$ is defined.
\end{defn}

In the situation of the above definition, we will say that the
extension $L/L_0$ is defined over the subextension $E/E_0$. If
$L/L_0$ is defined over $E/E_0$, then $E$ and $L_0$ are linearly
disjoint over $E_0$.

\begin{lemma} \label{l2.1}
Let $L/L_0$ be a finite separable extension of fields containing $K$.
If $L/L_0$ is defined over $E/E_0$, then $E$ is separable over $E_0$.
\end{lemma}

\begin{proof}
Recall that, in case $\fn{char}K=p>0$, $E$ is separable over $E_0$
if and only if $E^pE_0=E$ \cite[Corollary 7.6, p.293]{La}.

Let $E'=E^pE_0$. We will work in the category of vector spaces
over $E_0$. Define $V=E/E'$. From the short exact sequence $0\to
E'\to E\to V\to 0$, we get $0\to E'\otimes_{E_0} L_0\to
E\otimes_{E_0}L_0 \to V\otimes_{E_0}L_0 \to 0$.

Since $E$ and $L_0$ are linearly disjoint over $E_0$, it follows
that $E\otimes_{E_0}L_0\simeq EL_0$ and $E'\otimes_{E_0}L_0\simeq
E'L_0$. As $L$ is separable over $L_0$, we find that $L=L^pL_0$.
Thus $EL_0=L=L^pL_0=(EL_0)^pL_0=E^p(E_0L_0)=E'L_0$. Thus the map
$E'\otimes_{E_0}L_0\to E\otimes_{E_0} L_0$ is surjective. Hence
$V\otimes_{E_0} L_0=0$; it follows that $V=0$ because every
$E_0$-vector space is free over $E_0$ and therefore is faithfully
flat \cite[Chapter 1, Section 3]{Bo}.
\end{proof}

\begin{defn}\label{d2.2} \rm
Let $L$ be a field containing $K$ and $G$ be a finite group. $L$
is called a $G$-field over $K$ if there is a group homomorphism
$G\to \fn{Aut}_K(L)$, i.e. $G$ acts on $L$ by $K$-automorphisms;
since the ground field $K$ is understood throughout this article,
we simply call it a $G$-field. An intermediate field $E$,
$K\subset E\subset L$, is called a $G$-subfield of $L$ if $E$ is
invariant under the action of $G$, i.e. $\sigma(E)\subset E$ for
any $\sigma\in G$. A $G$-subfield $E$ is called a faithful
$G$-subfield if the restriction to $E$ of any non-identity element
of $G$ is not the identity map on $E$.
\end{defn}

With the aid of Lemma 2.2, the proof of the following two lemmas
are the same as that in \cite[Lemma 2.2 and Lemma 2.3]{BR}.

\begin{lemma}\label{l2.3}
Let $L$ be a field containing $K$,
and $G$ be a finite subgroup of $\fn{Aut}_K(L)$.

\r{(1)} If $E$ is a $G$-subfield of $L$, then $L/L^G$ is defined over $E/E^G$ if and only if $E$ is a faithful
$G$-subfield.

\r{(2)} There is always a faithful $G$-subfield $E$ of $L$ such that $\fn{ed}_K(L/L^G)=\fn{trdeg}_K E$.
\end{lemma}

\begin{lemma}\label{l2.4}
Let $L/L_0$ be a finite separable extension of fields containing $K$ and $L^\#$ be the normal closure of $L$ over $L_0$.
Then $\fn{ed}_K(L/L_0)=\fn{ed}_K(L^\#/L_0)$.
\end{lemma}

\begin{defn}\label{d2.5}\rm
Let $G$ be a finite group.
Define $\fn{ed}_K(G)=\fn{ed}_K(K(V)/K(V)^G)$ where $G\to GL(V)$ is a faithful representation of $G$ and
$K(V)$ is the function field of the affine space $V$ over $K$.
(Explicitly, if $x_1,\ldots,x_n$ is a basis of the dual space of $V$,
then $K(V)=K(x_1,\ldots,x_n)$ is the rational function field of $n$ variables over $K$.)
\end{defn}

Note that if $\rho:G\to GL(V_{reg})$ is the regular representation of the group $G$,
then $\rho$ is a faithful representation.

If $G\to GL(V)$ and $G\to GL(W)$ are two faithful representation
of a finite group $G$, we will show that
$\fn{ed}_K(K(V)/K(V)^G)=\fn{ed}_K(K(W)/K(W)^G)$; thus the
definition of $\fn{ed}_K(G)$ is independent of the choice of the
faithful representation of $G$. Before proving this fact, we
recall the following two theorems.

\begin{theorem}\label{t2.6}
\r{(Hajja and Kang \cite[Theorem 1]{HK})}
Let $G$ be a finite group acting on $L(x_1,\ldots,x_n)$,
the rational function field of $n$ variables over a field $L$.
Suppose that

\r{(i)} for any $\sigma\in G$, $\sigma(L)\subset L$;

\r{(ii)} the restriction of the actions of $G$ on $L$ is faithful;

\r{(iii)} for any $\sigma\in G$,
\[
\left(\begin{array}{c} \sigma(x_1) \\ \vdots \\ \sigma(x_n) \end{array}\right) =
A(\sigma)\left(\begin{array}{c} x_1 \\ \vdots \\ x_n \end{array}\right)+B(\sigma)
\]
where $A(\sigma)\in GL_n(L)$ and $B(\sigma)$ is an $n\times 1$ matrix over $L$.

Then there exist $z_1,\ldots,z_n\in L(x_1,\ldots,x_n)$ such that
$L(x_1,\ldots,x_n)=L(z_1,\ldots,z_n)$ with $\sigma(z_i)=z_i$ for
any $\sigma\in G$, any $1\le i \le n$.
\end{theorem}

\begin{theorem}\label{t2.7}
\r{(Jensen, Ledet and Yui \cite[Proposition 8.2.5, p.191]{JLY})}
Let $K$ be any field, $L$ be a field extension of $K$ and $L(x)$
be the rational function field of one variable over $L$. Let $E$
be an intermediate subfield with $K\subset E\subset L(x)$. Suppose
that $G$ is a finite group acting on $L(x)$ by $K$-automorphisms
such that \r{(i)} for any $\sigma\in G$, $\sigma(x)=x$, and
\r{(ii)} both $L$ and $E$ are $G$-subfield of $L(x)$. If
$\fn{trdeg}_KE\le \fn{trdeg}_KL<\infty$, then there is a
$K$-linear $G$-embedding of the field $E$ into $L$, i.e. there is
a $K$-morphism of fields $f:E\to L$ such that $f(\sigma\cdot
y)=\sigma\cdot f(y)$ for any $\sigma\in G$, any $y\in E$.
\end{theorem}

Note that Theorem \ref{t2.7} is an equivariant version (and a
generalization) of a theorem of Ohm \cite{Oh}.

It is time to show that
$\fn{ed}_K(K(V)/K(V)^G)=\fn{ed}_K(K(W)/K(W)^G)$ if $G\to GL(V)$
and $G\to GL(W)$ are faithful representations of $G$.

Let $x_1,\ldots,x_m$ and $y_1,\ldots,y_n$ are bases of the dual
spaces of $V$ and $W$ respectively. Then $K(V)=K(x_1,\ldots,x_m)$,
$K(W)=K(y_1,\ldots,y_n)$ and $K(V\oplus
W)=K(x_1,\ldots,x_m,y_1,\ldots,y_n)$. Without loss of generality,
we may assume that $\fn{ed}_K(K(V)$ $/K(V)^G)\le
\fn{ed}_K(K(W)/K(W)^G)$. Let $E$ be a faithful $G$-subfields of
$K(V)$ such that $\fn{trdeg}_KE=\fn{ed}_K(K(V)/K(V)^G)$. Consider
the action of $G$ on $K(V\oplus
W)=K(x_1,\ldots,x_m,y_1,\ldots,y_n)$. Apply Theorem \ref{t2.6} and
write $K(V\oplus W)=K(W)(z_1,\ldots,z_m)$ where $\sigma(z_i)=z_i$
for any $\sigma\in G$, any $1\le i\le m$. Note that $E\subset
K(V)\subset K(V\oplus W)=K(W)(z_1,\ldots,z_m)$. Apply Theorem
\ref{t2.7}. We may regard $E$ as a faithful $G$-subfield of
$K(W)$. Thus $\fn{trdeg}_K E\ge \fn{ed}_K(K(W)/K(W)^G)$. Hence the
result.

\begin{prop}\label{p2.8}
Let $G$ be a finite group and $L$ be a  faithfully $G$-field over
$K$. Then $\fn{ed}_K(L/L^G)\le \fn{ed}_K(G)$.
\end{prop}

\begin{proof}
Let $G\to GL(V)$ be a faithful representation of $G$. We may
choose $L$ and $K(V)$ so that they are free over $K$, i.e.
$\fn{trdeg}_K L\cdot K(V)=\fn{trdeg}_K L+\fn{trdeg}_K K(V)$. Then
$L\cdot K(V)=L(x_1,\ldots,x_n)$ where $x_1,\ldots,x_n$ is a basis
of the dual space of $V$. By Theorem \ref{t2.6} write
$L(x_1,\ldots,x_n)=L(z_1,\ldots,z_n)$ where $\sigma(z_i)=z_i$ for
any $\sigma\in G$, any $1\le i\le n$. Let $E$ be a faithful
$G$-subfield of $K(V)$ with $\fn{trdeg}_K E=\fn{ed}_K(G)$. If
$\fn{trdeg}_KE>\fn{trdeg}_K L$, then $\fn{ed}_K(L/L^G)\le
\fn{trdeg}_KL< \fn{trdeg}_K E=\fn{ed}_K(G)$ as we desire.
Otherwise, apply Theorem \ref{t2.7} to embed $E$ into $L$ because
$K\subset E\subset K(V)\subset L\cdot K(V)=L(z_1,\ldots,z_n)$.
Thus $\fn{trdeg}_KE\ge\fn{ed}_K(L/L^G)$ and we find that
$\fn{ed}_K(L/L^G)\le ed_K(G)$.
\end{proof}

The proof of the following lemma is the same as that in \cite{BR} and is omitted.

\begin{lemma}\label{l2.9}
Let $K$ be an arbitrary field and $G$ be a finite group.

\r{(1)} $\fn{ed}_K(G)=0$ if and only if $G=\{1\}$.

\r{(2)} If $G=G_1\times G_2$, then $\fn{ed}_K(G)\le \fn{ed}_K(G_1)+\fn{ed}_K(G_2)$.

\r{(3)} If $H$ is a subgroup of $G$, then $\fn{ed}_K(H) \le \fn{ed}_K(G)$.

\r{(4)} If $K'$ is a field extension of $K$, then $\fn{ed}_K(G)\ge \fn{ed}_{K'}(G)$.
\end{lemma}

\section{Fields of cross-ratios}

Let $n\ge 4$ and $S_n$ be the symmetric group acting on
$K(x_1,\ldots,x_n)$, the rational function field of $n$ variables
over $K$, by $\sigma\cdot x_i=x_{\sigma(i)}$ for any $\sigma\in
S_n$, any $1\le i\le n$.

If $i$, $j$, $k$, $l$ are four distinct elements in $\{1,2,\ldots,n\}$,
define the cross-ratio $[i,j;k,l]$ by
$$
[i,j;k,l]=(x_i-x_k)(x_j-x_l)(x_i-x_l)^{-1}(x_j-x_k)^{-1}\in K(x_1,\ldots,x_n).
$$

Let $L$ be the subfield of $K(x_1,\ldots,x_n)$ generated by
$[i,j;k,l]$ over $K$ where $i$, $j$, $k$, $l$ runs over all four
distinct integers chosen from the set $\{1,2,\ldots,n\}$. $L$ will
be called the field of cross-ratios over $K$.

\begin{lemma}\label{l3.1}
If $n\ge 5$ and $K$ is an arbitrary field, then the field of
cross-ratios $L$ is a faithful $S_n$-subfield of
$K(x_1,\ldots,x_n)$. Moreover, $L=K([1,2;3,i]:4\le i \le n)$ is
purely transcendental of degree $n-3$ over $K$.
\end{lemma}

\begin{proof}
Step 1.
For any $\sigma\in S_n$, $\sigma\cdot[i,j;k,l]=[\sigma(i),\sigma(j);\sigma(k),\sigma(l)]$.
Hence $L$ is a $S_n$-subfield of $K(x_1,\ldots,x_n)$.

If $n\ge 5$, then $S_n$ acts faithfully on $L$ because the
candidates of the kernel of the map $S_n\to \fn{Aut}_K(L)$ are
$S_n$, $A_n$ and $\{1\}$. It is easy to see that the first two
possibilities are not the kernel.

Step 2. Let $t=[1,2;3,4]$ and $\sigma\in S_n$ such that $\sigma$
leaves \{1,2,3,4\} invariant. Then $K(\sigma(t))=K(t)$. In fact,
it is routine to check that (i) $(12)\cdot t=(3,4)\cdot
t=\frac{1}{t}$, (ii) $(23)\cdot t=1-t$, (iii) $t$ is fixed by
(12)(34), (13)(24), (14)(23).

Step 3. For $4\le i \le n$, note that
$K(x_1,x_2,x_3,x_i)=K(x_1,x_2,x_3,[1,2;3,i])$. Thus
$K(x_1,x_2,\ldots,x_n)=K(x_1,x_2,x_3,[1,2;3,i]:4\le i \le n)$.
Define $L_0=K([1,2;3,i]:4\le i \le n)$. Then $L_0$ is purely
transcendental of degree $n-3$ over $K$ since $\fn{trdeg}_K
L_0(x_1,x_2,x_3)$ $=n$. It remains to show that $L_0=L$.

Step 4. If $i$, $j$, $k$, $l$ are any distinct four integers
chosen from $\{1,2,\ldots, n\}$, we will show that $[i,j;k,l]\in
L_0$.

Note that $[1,2;3,4][1,2;3,5]^{-1}=[1,2;5,4]$.
Using this formula, we will show that $[i,j;k,l]\in L_0$.

For example, if $\{i,j,k,l\}\cap \{1,2,3\}=\{1,2,3\}$, then we may
assume that $\{i,j,k,l\}=\{1,2,3,i\}$ without loss of generality.
Obviously $[1,2;3,i]\in L_0$ by definition. All other
cross-ratios, e.g. $[1,2;,i,3]$, $[1,i;2,3]$ etc, belong to $L_0$
by Step 2.

If $\{i,j,k,l\}\cap \{1,2,3\}=\{1,2\}$,
we may assume that $\{i,j,k,l\}=\{1,2,4,5\}$.
Thus $[1,2;4,5]=[1,2;3,5][1,2;3,4]^{-1}\in L_0$.

If $\{i,j,k,l\}\cap \{1,2,3\}=\{1\}$,
we may assume that $\{i,j,k,l\}=\{1,4,5,6\}$.
Then $[1,4;5,6]=[1,4;3,6][1,4;3,5]^{-1}$ while $[1,4;3,5]$, $[1,4;3,6]\in L_0$ by the previous case.

If $\{i,j,k,l\}\cap\{1,2,3\}=\emptyset$,
we may assume that $\{i,j,k,l\}=\{4,5,6,7\}$.
Then $[4,5;6,7]=[4,5;1,7][4,5;1,6]^{-1}$ and $[4,5;1,6]$, $[4,5;1,7]\in L_0$ as before.
\end{proof}

\begin{prop}\label{p3.2}
Let $K$ be an arbitrary field.
If $n\ge 5$, then $\fn{ed}_K(S_n)\le n-3$.
\end{prop}

\begin{proof}
Let $V=\bigoplus_{1\le i\le n} K\cdot x_i$ and $S_n$ act faithfully on $V$ by $\sigma\cdot x_i=x_{\sigma(i)}$ for any
$\sigma\in S_n$, any $1\le i \le n$.
Then $\fn{ed}_K(S_n)=\fn{ed}_K(K(x_1,\ldots,x_n)/K(x_1,\ldots,x_n)^{S_n})$.
Let $L$ be the field of cross-ratios over $K$.
Then $L\subset K(x_1,\ldots,x_n)$ and $L$ is a faithful $S_n$-subfield by Lemma \ref{l3.1}.
Since $\fn{trdeg}_K L=n-3$ by Lemma \ref{l3.1},
it follows that $\fn{ed}_K(S_n)\le n-3$.
\end{proof}

\begin{remark}\rm
This proposition was proved when $\fn{char}K=0$ in \cite[Theorem 6.5]{BR} and when $\fn{char}K\ne 2$ \cite[p.318--319]{BF}.
\end{remark}

\section{Central extensions}

In this section we will prove a generalization of \cite[Theorem
5.3]{BR}. As in the proof of Buhler and Reichstein, we will use
some basic results in valuation rings. The terminology of
valuations, valuation rings we adopt here follows those in
Bourbaki's book \cite[Chapter 6]{Bo}.

\begin{lemma}\label{l4.1}
Let $R$ be a rank-one discrete valuation ring $\rm {(}$DVR$\rm{)}$
containing a field $K$ and $\frak{M}$ be the maximal ideal of $R$.
If $\sigma\in\fn{Aut}_K(R)$ satisfying that \r{(i)}
$1<\fn{ord}(\sigma)<\infty$, and \r{(ii)} the induced
automorphisms of $\sigma$ on $R/\frak{M}$ and
$\frak{M}/\frak{M}^2$ are the identity maps, then $\fn{char}K=l>0$
and $\fn{ord}(\sigma)=l^r$ for some positive integer $r$.
\end{lemma}

\begin{proof}
The proof is adapted from \cite[Lemma 8.2.12]{JLY}.

Let $n=\fn{ord}(\sigma)$. If $\fn{char}K=0$, define an integer $m$
by $m=n$; if $\fn{char}K=l>0$, define the integer $m$ by $m=n'$
where $n=l^r\cdot n'$ with $l\nmid n'$. We will show that $m=1$.

Suppose not.
Define $\tau=\sigma^{n/m}$.
Then $\tau\in\fn{Aut}_K(R)$ and $\fn{ord}(\tau)=m>1$.
Choose $x\in R$ such that $\tau(x)\ne x$ and define $y=\tau(x)-x\ne 0$.

Since the induced actions of $\sigma$ on $R/\frak{M}$ and
$\frak{M}/\frak{M}^2$ are the identity maps, it is not difficult
to show that, for any $r=0,1,2,\ldots$, for any $u\in \frak{M}^r$,
it is necessary that $\sigma(u)-u\in \frak{M}^{r+1}$. The same is
valid for $\tau$, i.e. $\tau(u)-u\in\frak{M}^{r+1}$.

In particular, if $y\in\frak{M}^s\backslash\frak{M}^{s+1}$, then
$\tau(y)-y\in\frak{M}^{s+1}$.

From $\tau(x)=x+y$, we find that
$\tau^2(x)=\tau(x)+\tau(y)=x+y+y+v$ where $\tau(y)=y+v$ with
$v\in\frak{M}^{s+1}$. Hence $\tau^2(x)=x+2y$ (mod
$\frak{M}^{s+1}$). Proceed by induction. We find that
$\tau^i(x)=x+iy$ (mod $\frak{M}^{s+1}$) for $1\le i \le m$.
Consequently $\tau^m(x)=x+my$ (mod $\frak{M}^{s+1}$).

Since $\tau^m=1$, we find that $x=x+my$ (mod $\frak{M}^{s+1}$). It
follows that $my\in \frak{M}^{s+1}$. Since
$y\in\frak{M}^s\backslash\frak{M}^{s+1}$, it follows that $m\in
\frak{M}$. But $m$ is invertible in $K$ by assumption. Thus we get
a contradiction.
\end{proof}

Now let $R$, $\frak{M}$, $K$ be the same as in Lemma \ref{l4.1}.
Suppose that $G \subset \fn{Aut}_K(R)$ is a finite group. Define
\begin{eqnarray*}
G_0 &=& \{\sigma\in G:\mbox{ The induced automorphism of $\sigma$ on $R/\frak{M}$ is the identity map}\}, \\
G_1 &=& \{\sigma\in G:\mbox{ The induced automorphisms of $\sigma$ on $R/\frak{M}$ and $\frak{M}/\frak{M}^2$ are} \\
& & \hspace*{1cm}\mbox{ the identity maps}\}.
\end{eqnarray*}

Clearly $G_0$ and $G_1$ are normal subgroups. We will prove that,
under suitable conditions, $G_1=\{1\}$ and $G_0$ is contained in
the center of $G$. We will denote by $Z(G)$ the center of the
group $G$ in the sequel.

Since $\frak{M}$ is a principal ideal, we may write
$\frak{M}=\langle t\rangle$.

For any $\sigma\in G$, choose an element $\lambda_\sigma\in R$ so
that $\sigma(t)=\lambda_\sigma t$  (because
$\sigma(\frak{M})=\frak{M}$ for any $\sigma\in G$). Denote by
$\overline{\lambda}_\sigma$ the image of $\lambda_\sigma$ in
$R/\frak{M}$. We record some properties of
$\overline{\lambda}_\sigma$ in the following lemma.

\begin{lemma}\label{l4.2}
\r{(1)} For any $\sigma_1,\sigma_2\in G$,
$\lambda_{\sigma_1\sigma_2}=\lambda_{\sigma_1}\cdot\sigma_1(\lambda_{\sigma_2})
\, \rm{(}$mod $\frak{M} \rm{)}$. In particular,
$\lambda_{\tau\sigma\tau^{-1}}=\tau(\lambda_\sigma) \, \rm{(}$mod
$\frak{M} \rm{)}$ for any $\sigma\in G_0$, any $\tau\in G$.

\r{(2)} The map $\Phi:G_0\to (R/\frak{M})^{\times}$ defined by
$\Phi(\sigma)=\overline{\lambda}_\sigma$ is a group homomorphism
where $(R/\frak{M})^{\times}$ is the multiplicative group of the
field $R/\frak{M}$.
\end{lemma}

\begin{proof}
(1) Modulo $\frak{M}^2$,
we have $\lambda_{\sigma_1\sigma_2}t=\sigma_1\sigma_2(t)=\sigma_1(\sigma_2(t))=\sigma_1(\lambda_{\sigma_2}t)
=\sigma_1(\lambda_{\sigma_2})\cdot\sigma_1(t)=\sigma_1(\lambda_{\sigma_2})\cdot\lambda_{\sigma_1}t
=[\lambda_{\sigma_1}\cdot\sigma_1(\lambda_{\sigma_2})]t$.
Hence the result.

Suppose that $\sigma\in G_0$ and $\tau\in G$. Modulo $\frak{M}$,
we have
$\lambda_{\tau\sigma\tau^{-1}}=\lambda_\tau\cdot\tau(\lambda_{\sigma\tau^{-1}})
=\lambda_\tau\cdot
\tau(\lambda_\sigma\cdot\sigma(\lambda_{\tau^{-1}}))
=\tau(\lambda_\sigma)\cdot\lambda_\tau\cdot\tau\sigma(\lambda_{\tau^{-1}})=\tau(\lambda_\sigma)\cdot\lambda_\tau\cdot\tau(\lambda_{\tau^{-1}})$
because $\sigma$ is the identity map on $R/\frak{M}$. On the other
hand, note that $\lambda_1=1$. Thus
$1=\lambda_1=\lambda_{\tau\tau^{-1}}=\lambda_\tau\cdot\tau(\lambda_{\tau^{-1}})$
(mod $\frak{M}$). It follows that
$\lambda_{\tau\sigma\tau^{-1}}=\tau(\lambda_\sigma)$ (mod
$\frak{M}$).

(2) It suffices to show that, for any $\sigma_1,\sigma_2\in G_0$,
$\lambda_{\sigma_1\sigma_2}=\lambda_{\sigma_1}\cdot\lambda_{\sigma_2}$
(mod $\frak{M}$). By (1), we have
$\lambda_{\sigma_1\sigma_2}=\lambda_{\sigma_1}\cdot\sigma_1(\lambda_{\sigma_2})
=\lambda_{\sigma_1}\cdot\lambda_{\sigma_2}$ (mod $\frak{M}$) since
$\sigma_1$ is the identity map on $R/\frak{M}$.
\end{proof}

\begin{lemma}\label{l4.3}
Let $R$ be a DVR containing a field $K$, $\frak{M}$ be its maximal ideal.
Suppose that $G$ is a finite subgroup of $\fn{Aut}_K(R)$ and define $G_0$ and $G_1$ as above.
Assume furthermore that $K$ is algebraically closed in the field $R/\frak{M}$.
Then

\r{(i)} $\tau\sigma\tau^{-1}\sigma^{-1}\in G_1$ for any $\tau\in G$, any $\sigma\in G_0$;

\r{(ii)} $G_0/G_1$ is a cyclic group;

\r{(iii)} if the order of $G_0/G_1$ is $m$, then $\zeta_m\in K$; in particular, $\fn{char}K\nmid m$.
\end{lemma}

\begin{coro}\label{c4.4}
Let the assumptions be the same as in Lemma \ref{l4.3}. Assume
that \r{(i)} $\fn{char}K=0$, or \r{(ii)} $\fn{char}K=l>0$ and $G$
contains no non-trivial normal $l$-subgroup. Then $G_1=\{1\}$,
$G_0$ is a cyclic group and $G_0\subset Z(G)$.
\end{coro}

\renewcommand{\proofname}{\indent Proof of Corollary \ref{c4.4}}
\begin{proof}
Assume that Lemma \ref{l4.3} is valid.

Note that $G_1=\{1\}$. Otherwise, $\fn{char}K=l>0$ by Lemma
\ref{l4.1} and $G_1$ is a \mbox{$l$-subgroup} of $G$. Note that
$G_1$ is normal in $G$. This leads to a contradiction to our
assumption.

Now apply Lemma \ref{l4.3}.
Since $G_1=\{1\}$, we find $\sigma\tau=\tau\sigma$ for any $\sigma\in G_0$, any $\tau\in G$.
Thus $G_0\subset Z(G)$.
Again by Lemma \ref{l4.3} we know that $G_0\simeq G_0/G_1$ is a cyclic group.
\end{proof}

\renewcommand{\proofname}{\indent Proof of Lemma \ref{l4.3}}
\begin{proof} ~

Consider the group homomorphism $\Phi:G_0\to
(R/\frak{M})^{\times}$ in Lemma \ref{l4.2}. $G_1$ is just the
kernel of $\Phi$. The image of $G_0$ is a cyclic group, say, of
order $m$. Let $\alpha\in R$ such that $\overline{\alpha}$ is a
generator of the image of $\Phi$.

Since $K\subset R/\frak{M}$ and $\overline{\alpha}\in R/\frak{M}$ is algebraic over $K$,
it follows that $\alpha\in K$ because we assume that $K$ is algebraically closed in $R/\frak{M}$.
Clearly $\fn{char}K\nmid m$ and $G$ is a cyclic group.

It remains to show that $\tau\sigma\tau^{-1}\sigma^{-1}\in G_1$ for any $\tau\in G$, $\sigma\in G_0$.
It suffices to show that $\lambda_{\tau\sigma\tau^{-1}\sigma^{-1}}=1$ (mod $\frak{M}$).

Since $\fn{Im}(\Phi)=\langle\alpha\rangle$, it follows that
$\overline{\lambda}_\sigma=\alpha^i\in K$ for some $i$. Since
$\tau\sigma\tau^{-1}\in G_0$, it follows that
$\overline{\lambda}_{\tau\sigma\tau^{-1}}=\tau(\overline{\lambda}_\sigma)$
by Lemma \ref{l4.2}. Thus
$\tau(\overline{\lambda}_\sigma)=\tau(\alpha^i)=\alpha^i$ because
$G$ is a $K$-linear automorphism group. Thus
$\overline{\lambda}_{\tau\sigma\tau^{-1}\sigma^{-1}}
=\overline{\lambda}_{\tau\sigma\tau^{-1}}\cdot\overline{\lambda}_{\sigma^{-1}}=\alpha^i\cdot
\alpha^{-i}=1$.
\end{proof}

\begin{theorem}\label{t4.5}
Let $p$ be a prime number, $G$ be a finite group.
Assume that the following four properties are valid,

\r{(i)}
$K$ is a field satisfying that either $\fn{char}K=0$ or $\fn{char}K=l>0$ such that
$G$ contains no non-trivial normal $l$-subgroup,

\r{(ii)}
there is an element $\sigma\in Z(G)\subset G$ with $\fn{ord}(\sigma)=p$,

\r{(iii)} there is a linear character $\chi: G\to K^{\times}$ with
$\chi(\sigma)\ne 1$,

\r{(iv)} for any $\tau\in Z(G)$ with $\fn{ord}(\tau)=m$, it is
necessary that $\zeta_m\notin K$ if $\langle\sigma\rangle$ is a
proper subgroup of $\langle\tau\rangle$.

Then $\fn{ed}_K(G)=\fn{ed}_K(G/\langle\sigma\rangle)+1$.
\end{theorem}

We will consider also a split version of Theorem \ref{t4.5}, i.e.

\begin{theorem}\label{t4.6}
Let $p$ be a prime number, $G=G'\times\langle\sigma\rangle$ be a finite group.
Assume that the following four properties are valid.

\r{(i)} $K$ is a field satisfying that either $\fn{char}K=0$ or
$\fn{char}K=l>0$ such that $G$ contains no non-trivial normal
$l$-subgroup,

\r{(ii)} $\fn{ord}(\sigma)=p$,

\r{(iii)} $\zeta_p\in K$,

\r{(iv)} for any prime factor $p'$ of $|Z(G')|$ such that $p' \neq
p$, it is necessary that $\zeta_{p'}\notin K$.

Then $\fn{ed}_K(G)=\fn{ed}_K(G')+1$.
\end{theorem}

\begin{remark} \rm
In Theorem \ref{t4.5}(ii), since $\chi(\sigma)\ne 1$ and $\fn{ord}(\sigma)=p$,
it follows that $\chi(\sigma)$ is a primitive $p$-th root of unity.
In particular, $\zeta_p\in K$ and $\fn{char} K\nmid p$.
\end{remark}

\renewcommand{\proofname}{\indent Proof of Theorem \ref{t4.5} and Theorem \ref{t4.6}}
\begin{proof}~

In the situation of Theorem \ref{t4.5} write
$G'=G/\langle\sigma\rangle$. In the situation of Theorem
\ref{t4.6} define a linear character $\chi: G=G'\times
\langle\sigma\rangle \to K^{\times}$ by
$\chi(\tau\sigma^i)=\zeta^i_p$ for any $\tau\in G'$, any $0\le
i\le p-1$. Thus there is an embedding of $G$ into $G'\times \chi
(G)$ for both situations. By Lemma \ref{l2.9}, $\fn{ed}_K(G)\le
\fn{ed}_K(G'\times \chi(G))\le \fn{ed}_K(G')+1$ because
$\fn{ed}_K(\chi(G))=1$. It remains to show that $\fn{ed}_K(G)\ge
\fn{ed}_K(G')+1$.

Step 1. Let $G'\to GL(V)$ be a faithful representation of $G'$.
Then $G$ acts faithfully on $V\oplus K$ by $\tau\cdot
(v+\lambda)=\tau\cdot v+\chi(\tau)\lambda$ for any $\tau \in G$,
any $v\in V$, any $\lambda\in K$.

Write the function field $K(V\oplus K)=K(V)(x)$.
Thus $\tau\cdot x=\chi(\tau)^{-1}\cdot x$.
Hence $\tau\cdot(\frac{1}{x})=\chi(\tau)\cdot(\frac{1}{x})$.
By abusing the notation, we will write $K(V\oplus K)=K(V)(x)$ with $\tau\cdot x=\chi(\tau)\cdot x$.

Step 2.
Let $E$ be a faithful $G$-subfield of $K(V)(x)$ with $\fn{ed}_K(G)=\fn{trdeg}_KE$.
We will prove that $\fn{trdeg}_KE\ge \fn{ed}_K(G')+1$.

Let $\nu$ be a rank-one discrete valuation on $K(V)(x)$ determined by the prime ideal $\langle x\rangle$
in the polynomial ring $K(V)[x]$.
Let $O_\nu$ be its valuation ring and $\frak{M}_\nu$ be its maximal ideal.
Then the field $K(V)$ maps isomorphically onto the residue field of $\nu$,
i.e. $K(V)\simeq O_\nu/\frak{M}_\nu$.
In particular, $K$ is algebraically closed in $O_\nu/\frak{M}_\nu$.

Note that $O_\nu$ is invariant under the action of $G$. In fact,
$K[V]$ is isomorphic to a polynomial ring and $K[V]$ is invariant
under the action of $G$; moreover, every element of $O_\nu$ can be
written as
$u=x^n(a_0+a_1x+\cdots+a_rx^r)(b_0+b_1x+\cdots+b_sx^s)^{-1}$ where
$n\ge 0$, $a_i,b_j\in K[V]$ with $a_0b_0\ne 0$. Thus for any
$\tau\in G$, $\tau\cdot u=(\chi(\tau)x)^n(\tau(a_0)+$
$\chi(\tau)\cdot\tau(a_1)\cdot
x+\cdots+\chi(\tau)^r\cdot\tau(a_r)\cdot x^r)
(\tau(b_0)+\chi(\tau)\cdot\tau(b_1)\cdot
x+\cdots+\chi(\tau)^s\cdot\tau(b_s)\cdot x^s)^{-1}\in O_\nu$.

Step 3.
Consider the restriction of $\nu$ on the field $E$.

We claim that $\nu$ is non-trivial on $E$. Otherwise, $E\subset
O_\nu$. Thus there is a $G$-embedding of $E$ into $K(V)$, because
$K(V)$ is isomorphic to the residue field of $\nu$. However, the
element $\sigma\in G$ acts trivially on $K(V)$ while $G$ is
faithful on $E$. This is a contradiction.

Let $R$ be the valuation ring of the restriction of $\nu$ on $E$ and $\frak{M}$ be its maximal ideal.
Write $\bm{E}=R/\frak{M}$,
the residue field.
Since $R=E\cap O_\nu$,
it follows that $R$ is invariant under the action of $G$.
Hence there is a natural $G$-embedding of $\bm{E}$ into $K(V)$.
Note that $\fn{trdeg}_K\bm{E}\le\fn{trdeg}_KE-1$ by \cite[p.439, Corollary 1]{Bo} or by \cite[Lemma 8.1.2, p.187]{JLY}.

The action of $\sigma$ on $V$ is trivial while $G'=G/\langle\sigma\rangle$ acts faithfully on $V$.
We will show that $G'$ acts faithfully also on $\bm{E}$.
Thus $\bm{E}$ is a faithful $G'$-subfield of $K(V)$.
It follows that $\fn{ed}_K(G')\le\fn{trdeg}_K(\bm{E})\le\fn{trdeg}_K E-1$,
which will conclude our proof of Theorem \ref{t4.5} and Theorem \ref{t4.6}.

Step 4.
We will prove that $G'$ acts faithfully on $\bm{E}$.

Define $G_0=\{\tau\in G:$ The induced action of $\tau$ on
$R/\frak{M}\simeq \bm{E}$ is the identity map$\}$. Clearly
$\sigma\in G_0$. Note that $G'$ acts faithfully on $\bm{E}$ if and
only if $G_0=\langle\sigma\rangle$.

Since $E$ is a faithful $G$-subfield, the action of $G$ on $R$ is
faithful. Thus $G \subset \fn{Aut}_K(R)$ and $K$ is algebraically
closed in $\bm{E}$ because so is it in $K(V)$ which contains
$\bm{E}$.

Define $G_1=\{\tau\in G:$ The induced actions of $\tau$ on
$R/\frak{M}$ and $\frak{M}/\frak{M}^2$ are the identity maps$\}$.
By Corollary \ref{c4.4}, $G_1=\{1\}$ and $G_0$ is a cyclic
subgroup of $Z(G)$. Let $m=|G_0|$. Then $\zeta_m\in K$ by Lemma
\ref{l4.3}.

We will show that $G_0=\langle\sigma\rangle$, which is equivalent
to showing that $p=m$. Note that $m$ is divisible by $p$.

In the situation of Theorem \ref{t4.5},
if $p<m$, the property $\zeta_m\in K$ violates the assumption (iv) in Theorem \ref{t4.5}.

In the situation of Theorem \ref{t4.6}, since $G_0$ is cyclic
containing $\sigma$, it follows that
$G_0=H\times\langle\sigma\rangle$ for some subgroup $H\subset
Z(G')$. Should $p$ divide $|H|$, then $G_0$ would contain an
elementary abelian group of order $p^2$. But $G_0$ is a cyclic
group! Thus $p\nmid |H|$. Any prime divisor of $|H|$ is a divisor
of $\gcd\{m,|Z(G')|\}$. If $p<m$, the property $\zeta_m\in K$
would lead to a contradiction to the assumption (iv) of Theorem
\ref{t4.6}.
\end{proof}

\begin{theorem}\label{t4.7}
If $p$ is a prime number and $K$ is a field such that $\zeta_p\in
K$, then $\fn{ed}_K((\bm{Z}/p\bm{Z})^r)=r$. In particular, for a
field $K$ with $\fn{char}K\ne 2$,
$\fn{ed}_K((\bm{Z}/2\bm{Z})^r)=r$.
\end{theorem}

\renewcommand{\proofname}{\indent Proof}
\begin{proof}
Apply Theorem \ref{t4.5} or Theorem \ref{t4.6}.
\end{proof}

\section{Symmetric groups}

We will study the essential dimensions of symmetric groups $S_n$
and alternating groups $A_n$ in this section.

First recall a general version of L\"{u}roth's Theorem.

\begin{theorem}\label{t5.1}
Let $K$ be an arbitrary field and $K(x_1,\ldots,x_n)$ be
the rational function field of $n$ variables over $K$.
If $E$ is an intermediate field $K\subset E\subset K(x_1,\ldots,x_n)$
with $\fn{trdeg}_KE=1$,
then $E=K(y)$ for some element $y\in E$.
\end{theorem}

\begin{proof}
See \cite[p.87]{Na}. Alternatively we may apply Theorem \ref{t2.7}
by taking the trivial $G$-action on $L(x)$, i.e.\ the original
version  of Ohm's Theorem \cite{Oh}. (Note that, in Theorem
\ref{t2.7}, we don't require whether the action of $G$ on $L(x)$
is faithful or not.) Thus we may embed $E$ successively into
$K(x_1,\ldots,x_{n-1}), K(x_1,\ldots,x_{n-2}), \ldots, K(x_1)$.
Then apply the standard form of L\"{u}roth's Theorem.
\end{proof}

\begin{lemma}\label{l5.2}
Let $K$ be a field with $\fn{char}K=l>0$. Suppose that $\sigma\in
PGL_2(K)$ and $\fn{ord}(\sigma)$ in $PGL_2(K)$ is finite. Then
either $l\nmid \fn{ord}(\sigma)$ or $\fn{ord}(\sigma)=l$.
\end{lemma}

\begin{proof}
Note that the order of $\sigma$ in $PGL_2(K)$ is the same as that
in $PGL_2(\overline{K})$ where $\overline{K}$ is the algebraic
closure of $K$.

Choose a matrix $T\in GL_2(K)$ such that its image in $PGL_2(K)$
is $\sigma$. Find the Jordan canonical form of $T$ in
$GL_2(\overline{K})$. It is of the form
\[
\left(\begin{array}{cc}  a & 1 \\ 0 & a  \end{array}\right) \mbox{ \ or \ }
\left(\begin{array}{cc}  a & 0 \\ 0 & b  \end{array}\right)
\]
where $a,b\in\overline{K}\backslash\{0\}$. It is not difficult to
see that the order of the image of the above matrix in
$PGL_2(\overline{K})$ is $l$ for the first case, and the order of
the image in $PGL_2(\overline{K})$ is relatively prime to $l$ for
the second case.
\end{proof}

\begin{lemma}\label{l5.3}
\r{(\cite[Lemma 7.2]{BF})}
Let $K$ be an arbitrary field and $G$ be a finite group.
If $\fn{ed}_K(G)=1$, then $G$ may be embedded into $PGL_2(K)$.
\end{lemma}

\begin{proof}
Let $G\to GL(V)$ be a faithful representation.
If $\fn{ed}_K(G)=1$, there is a faithful $G$-subfield $E$ in $K(V)$ with $\fn{trdeg}_KE=1$.
By Theorem \ref{t5.1} $E=K(x)$ for some $x\in E$.
Thus $G$ may be regarded as a subgroup of $\fn{Aut}_K(K(x))\simeq PGL_2(K)$.
\end{proof}

\renewcommand{\proofname}{\indent Proof of Theorem \ref{t1.2}}
\begin{proof}~

Step 1.
Part (1) is just the same as Proposition \ref{p3.2}.

For the proof of Part (3), note that $\fn{ed}_K(S_6)\le 3$ by Part (1).
When $\fn{char}K\ne 2$, we find $\fn{ed}_K((\bm{Z}/2\bm{Z})^3)=3$ by Theorem \ref{t4.7}.
Since we may embed $(\bm{Z}/2\bm{Z})^3$ into $S_6$,
$3=\fn{ed}_K((\bm{Z}/2\bm{Z})^3)\le \fn{ed}_K(S_6)\le 3$.
Hence the result.

Step 2.
We will prove $\fn{ed}_K(S_4)=\fn{ed}_K(S_5)=2$.

By Part (1), $\fn{ed}_K(S_5)\le 2$.
Thus $\fn{ed}_K(S_4)\le\fn{ed}_K(S_5)\le 2$.

If $\fn{char}K\ne 2$,
we find that $\fn{ed}_K((\bm{Z}/2\bm{Z})^2)=2$ by Theorem \ref{t4.7}.
Hence $2=\fn{ed}_K((\bm{Z}/2\bm{Z})^2)\le \fn{ed}_K(S_4)\le \fn{ed}_K(S_5)\le 2$.
Thus $\fn{ed}_K(S_4)=\fn{ed}_K(S_5)=2$.

If $\fn{char}K=2$ and $\fn{ed}_K(S_4)=1$,
then $S_4$ can be embedded into $PGL_2(K)$ by Lemma \ref{l5.3}.
Since $S_4$ has an element of order 4,
it follows that $PGL_2(K)$ has an element of order 4,
which is impossible by Lemma \ref{l5.2}.
Thus $2\le\fn{ed}_K(S_4)\le \fn{ed}_K(S_5)\le 2$.

Step 3. We will prove $\fn{ed}_K(S_2)=\fn{ed}_K(S_3)=1$. We will
prove it by using \cite[Corollary 4.2]{BR}, i.e.\ the reduction to
the essential dimension of the general equation (see the next
paragraph).

Let $K$ be an arbitrary field, $t_1,t_2,\ldots,t_n$ be
algebraically independent over $K$, $L_0=K(t_1,\ldots,t_n)$. Let
$f(X)=X^n+t_1X^{n-1}+t_2X^{n-2}+\cdots+t_{n-1}X+t_n\in L_0[X]$ be
the general polynomial of degree $n$ over $K$. Let $\alpha$ be a
root of $f(X)=0$ and $L=L_0(\alpha)$. As in \cite[Corollary
4.2]{BR} it is easy to see that $\fn{ed}_K(L/L_0)=\fn{ed}_K(S_n)$.

If $\fn{char}K\nmid n$, then we may "kill" the coefficient of
$X^{n-1}$ for an equation of degree $n$ by the classical trick. So
we get an equation $g(X)=X^n+a_2X^{n-2}+a_3X^{n-3}+\cdots+a_n\in
L_0[X]$. We may re-scale the roots of $g(X)=0$ as in
\cite[p.195]{JLY}. Explicitly, substitute $\lambda X$ for $X$ in
$g(X)=0$ where the value of $\lambda$ is to be determined. We get
$(\lambda X)^n+a_2(\lambda X)^{n-2}+\cdots+a_{n-1}(\lambda
X)+a_n=0$, i.e.
$X^n+a_2\lambda^{-2}X^{n-2}+\cdots+a_{n-1}\lambda^{-(n-1)}X+a_{n-1}\lambda^{-n}=0$.
We get $\lambda=a_n/a_{n-1}$ by requiring
$a_{n-1}\lambda^{-(n-1)}=a_n\lambda^{-n}$. Hence we get an
equation
$h(X)=X^n+b_2X^{n-2}+b_3X^{n-3}+\cdots+b_{n-1}X+b_{n-1}=0$, i.e.\
only $n-2$ variables $b_2,b_3,\ldots,b_{n-1}$ are involved.

When $n \le 3$, the above process will not work for two
situations: (i) $n=2$ and $\fn{char}K=2$; (ii) $n=3$ and
$\fn{char}K=3$.

Consider the first case $n=2$ and $\fn{char}K=2$. We may re-scale
the roots of $f(X)=X^2+t_1X+t_2=0$. So $\fn{ed}_K(L/L_0)=1$ in
this case.

Consider the case of $n=3$ and $\fn{char}K=3$. We will "kill" the
coefficient of $X$ instead of that of $X^2$. Substitute
$X+\lambda$ for $X$ in $f(X)=X^3+t_1X^2+t_2X+t_3=0$. We get
$X^3+t_1X^2+(t_2-\lambda t_1)X+(t_3+\lambda t_2+\lambda^2 t_1)=0$.
Define $\lambda=t_2t^{-1}_1$ by requiring $t_2-\lambda t_1=0$. We
get the polynomial $g(X)=X^3+a_1 X^2+a_2\in L_0[X]$. If $\alpha$
is a root of $g(X)=0$, consider the minimum polynomial of
$\alpha^{-1}$. It is of the form $h(X)=X^3+b_1X+b_2\in L_0[X]$. We
may re-scale the roots of $h(X)=0$ and get an equation with only
one variable in the coefficients. Thus $\fn{ed}_K(L/L_0)=1$.
\end{proof}

\begin{theorem}\label{t5.4}
\r{(1)} Let $K$ be any field with $\fn{char}K\ne 2$.
Then
\begin{eqnarray*}
\r{ed}_K(S_{n+2}) &\ge& \r{ed}_K(S_n)+1 \mbox{ \ for any }n\ge 1,\mbox{ and} \\
\r{ed}_K(S_n) &\ge& \left\lfloor\frac{n}{2}\right\rfloor \mbox{ \ for any }n\ge 1.
\end{eqnarray*}

\r{(2)} If $\fn{char}K=2$ and $K\supset\bm{F}_4$, then
$\fn{ed}_K(S_{n+3})\ge \fn{ed}_K(S_n)+1$ for any $n\ne 4$.

\r{(3)} If $K$ is a field with $\fn{char}K=2$, then
$\fn{ed}_K(S_n)\ge \lfloor\frac{n+1}{3}\rfloor$ for $n\ge 1$.
\end{theorem}

\renewcommand{\proofname}{\indent Proof}
\begin{proof}
(1) Note that $S_n\times \bm{Z}/2\bm{Z}\hookrightarrow S_{n+2}$.
Thus we may apply Theorem \ref{t4.6} when $n\ge 4$. (Note that
$Z(S_n)=\{1\}$ when $n\ge 3$.) Hence $\fn{ed}_K(S_{n+2})\ge
\fn{ed}_K(S_n)+1$ for any $n\ge 4$. The situation when $n=1$, 2 or
3 is easy by Theorem \ref{t1.2}.

It follows that $\fn{ed}_K(S_n)\ge \lfloor\frac{n}{2}\rfloor$.
Done.

(2) Suppose that $K\supset\bm{F}_4$, i.e. $\zeta_3\in K$. Apply
Theorem \ref{t4.6} with $S_n\times \bm{Z}/3\bm{Z}\hookrightarrow
S_{n+3}$. Hence $\fn{ed}_K(S_{n+3})\ge \fn{ed}_K(S_n)+1$ for $n\ge
5$. The other cases (with $n\ne 4$) can be verified easily.

(3) Consider the case $K\supset\bm{F}_4$ first. Using (2) we may
prove that $\fn{ed}_K(S_n)\ge \lfloor\frac{(n+1)}{3}\rfloor$ for
$n\ge 8$. The situation when $n\le 7$ can be verified by using
Part (2) and Theorem 1.2.

If $K$ is any field with $\fn{char}K=2$, then $\fn{ed}_K(S_n)\ge
\fn{ed}_{K(\zeta_3)}(S_n)\ge \lfloor\frac{n+1}{3}\rfloor$.
\end{proof}

\begin{lemma}\label{l5.5}
\r{(1)} If $K$ is any field with $\fn{char}K=2$, then
$\fn{ed}_K(A_8)\le 3$.

\r{(2)} If $K$ is a field containing $\bm{F}_4$,
then $\fn{ed}_K(A_5)=1$.
\end{lemma}

\begin{remark} \rm
It will be shown in \cite{CHKZ} that, if $K$ is a field with
$\fn{char}K=2$ and $K$ doesn't contain $\bm{F}_4$, then
$\fn{ed}_K(A_4)=\fn{ed}_K(A_5)=2$.
\end{remark}

\begin{proof}
(1) Since $A_8\simeq GL_4(\bm{F}_2)$ (see \cite[p.363]{Ar} for example),
$A_8$ has a 4-dimensional faithful representation over $\bm{F}_2$ (and hence over $K$).
Thus $A_8$ acts faithfully on $K(x_1$, $x_2,x_3,x_4)$.
Define $E=K(y_1,y_2,y_3)$ where $y_i=x_i/x_4$ for $1\le i \le 3$.
Clearly $E$ is an $A_8$-subfield of $K(x_1,x_2,x_3,x_4)$.
Since $A_8$ is a simple group and the action of it on $E$ is non-trivial,
$A_8$ acts faithfully on $E$.
Thus $\fn{ed}_K(A_8)\le \fn{trdeg}_K E=3$.

(2) Let $K\supset\bm{F}_4$. Note that $A_5\simeq SL_2(\bm{F}_4)$
(see \cite[p.363]{Ar}). Thus we find that $\fn{ed}_K(A_5)\le 1$ as
in (1).
\end{proof}

\begin{theorem}\label{t5.6}
\r{(1)} Let $K$ be any field with $\fn{char}K\ne 2$.
Then
\begin{eqnarray*}
&&\r{ed}_K(A_3)=1,\ \r{ed}_K(A_4)=\r{ed}_K(A_5)=2, \\
&&\r{ed}_K(A_{n+4})\ge \r{ed}_K(A_n)+2 \mbox{ \ for }n\ge 4, \\
&&\r{ed}_K(A_n)\ge 2\left\lfloor\frac{n}{4}\right\rfloor \mbox{ \
for }n\ge 3.
\end{eqnarray*}

\r{(2)} Let $K$ be any field such that $K\supset\bm{F}_4$, then
$\fn{ed}_K(A_{n+3})\ge\fn{ed}_K(A_n)+1$ for $n\neq 4$.

\r{(3)} If $K$ is a field with $\fn{char}K=2$,
then $\fn{ed}_K(A_n)\ge \lfloor\frac{n}{3}\rfloor$ for $n\ge 3$.
\end{theorem}

\begin{proof}
(1) Since $\fn{ed}_K(S_3)=1$,
it follows that $\fn{ed}_K(A_3)=1$.
Similarly $\fn{ed}_K(A_5)\le 2$.
On the other hand,
$V_4\subset A_4$ where $V_4\simeq (\bm{Z}/2\bm{Z})^2$ is the Klein-four group.
Hence $2=\fn{ed}_K(V_4)\le \fn{ed}_K(A_4)\le \fn{ed}_K(A_5)\le 2$.

For $n\ge 5$, $A_n\times V_4\subset A_{n+4}$. Hence
$\fn{ed}_K(A_n\times V_4)\le \fn{ed}_K(A_{n+4})$. Write $A_n\times
V_4=G'\times \langle\sigma\rangle$ where
$G'=A_n\times\bm{Z}/2\bm{Z}$ and $\fn{ord}(\sigma)=2$. Apply
Theorem \ref{t4.5}. Note that $Z(G'\times \langle\sigma\rangle)$
has no cyclic subgroup containing $\langle\sigma\rangle$ properly.
Hence $\fn{ed}_K(A_n\times V_4)=\fn{ed}_K(G')+1$ by Theorem
\ref{t4.6}. On the other hand, $\fn{ed}_K(G')=\fn{ed}_K(A_n)+1$.
Hence $\fn{ed}_K(A_{n+4})\ge\fn{ed}_K(A_n)+2$ for $n\ge 5$. The
relation $\fn{ed}(A_8)\ge \fn{ed}(A_4)+2$ can be verified because
$V_4\times V_4\subset A_8$ and $\fn{ed}_K(V_4\times V_4)=4$.

The property $\fn{ed}_K(A_n)\ge 2\lfloor\frac{n}{4}\rfloor$ for
$n\ge 4$ can be proved by using the relation
$\fn{ed}_K(A_{n+4})\ge\fn{ed}_K(A_n)+2$ and by induction. The case
for $\fn{ed}_K(A_3)$ is obvious.

(2) Consider the field $K\supset \bm{F}_4$.

Note that $A_n\times \bm{Z}/3\bm{Z}\subset A_{n+3}$. If $n\ge 5$,
apply Theorem \ref{t4.6} to get $\fn{ed}_K(A_{n+3})\ge
\fn{ed}_K(A_n)+1$ for $n\ge 5$. If $n=3$, note that
$\fn{ed}_K(A_6)\ge \fn{ed}_K((\bm{Z}/3\bm{Z})^2)=2$ by Theorem
4.7. Hence the result.

(3) As in the proof of Theorem \ref{t5.4}(3),
it suffices to prove the formula when $K$ is a field containing $\bm{F}_4$.

By Lemma \ref{l5.5},
$\fn{ed}_K(A_5)=1$ for such a field $K$.
Hence $\fn{ed}_K(A_3)=\fn{ed}_K(A_4)=1$.
Thus $\fn{ed}_K(A_n)=\lfloor\frac{n}{3}\rfloor$ for $3\le n\le 5$.

Note that that $\fn{ed}_K(A_6)\ge 2$ and $\fn{ed}_K(A_9)\ge
\fn{ed}_K(A_6)+1\ge 3$ by Part (2). Now it is clear that
$\fn{ed}_K(A_n)\ge \lfloor\frac{n}{3}\rfloor$ for $6\le n\le 9$.

When $n\ge 10$, we may use the induction hypothesis and the result
in (2).
\end{proof}

Finally we will determine the field $K$ when $\fn{ed}_K(D_n)=1$.
Note that various cases of $\fn{ed}_K(G)=1$ when
$G=\bm{Z}/n\bm{Z}$ or $D_n$ were investigated in \cite[Section
7]{BF}.

\begin{lemma}\label{l5.7}
Let $K$ be an arbitrary field, $\sigma\in PGL_2(K)$ an element of
finite order. If $n=\fn{ord}(\sigma)$ and $\fn{char}K \nmid n$,
then $\zeta_n+\zeta^{-1}_n\in K$.
\end{lemma}

\begin{remark}\rm
The above lemma is a generalization of \cite[Lemma 7.7]{BF} where
it was required that $n$ is a prime number.
\end{remark}

\begin{proof}
Choose a matrix $T\in GL_2(K)$ such that its image in $PGL_2(K)$
is $\sigma$. Then the Jordan canonical form of $T$ is
\[
\left(\begin{array}{cc}  \lambda \zeta_n & 0 \\ 0 & \lambda
\end{array} \right)
\]
for some $\lambda$ in the algebraic closure of $K$.

Note that the rational canonical form of $T$ is
\[
\left(\begin{array}{cc}  \lambda \zeta_n & 0 \\ 0 & \lambda
\end{array} \right) \mbox{ \ or \ } \left(\begin{array}{cc}  0 & a
\\ 1 & b  \end{array} \right)
\]
where $\lambda,a,b\in K$ (according to whether the characteristic
polynomial is reducible over $K$ or irreducible over $K$). Then
first possibility will imply $\zeta_n\in K$; in particular
$\zeta_n+\zeta^{-1}_n\in K$. It remains to consider the second
possibility.

Compare the traces and the determinants of these two canonical
forms. We find that $b=\lambda\zeta_n+\lambda$ and
$-a=\lambda^2\zeta_n$. Thus
$\zeta_n+\zeta^{-1}_n+2=(\lambda\zeta_n+\lambda)^2(\lambda^2\zeta_n)^{-1}=b^2\cdot
(-a)^{-1}\in K$.
\end{proof}

\begin{theorem}\label{t5.8}
\r{(1)} Let $K$ be a field with $\fn{char}K\ne 2$. Then
$\fn{ed}_K(D_n)=1$ if and only if both the following conditions
are valid: \r{(i)} $n$ is an odd integer, and \r{(ii)}
$\zeta_n+\zeta^{-1}_n\in K$ if $\fn{char}K\nmid n$, or $n=p$ if
$\fn{char}K=p>0$ with $p \, | \, n$.

\r{(2)} Let $K$ be a field with $\fn{char}K=2$. Then
$\fn{ed}_K(D_n)=1$ if and only if $\zeta_n+\zeta^{-1}_n\in K$ if
$n$ is an odd integer, or $|K|\ge 4$ with $n=2$ if $n$ is an even
integer.
\end{theorem}

\begin{proof}
(1) Assume that $\fn{char}K\ne 2$ and $\fn{ed}_K(D_n)=1$.

By Lemma \ref{l5.3} we may embed $D_n$ into $PGL_2(K)$. Then
$PGL_2(K)$ contains an element $\sigma$ of order $n$. By Lemma
\ref{l5.2} we find that either $\fn{char}K\nmid n$ or
$\fn{char}K=p>0$ with $n=p$. In case $\fn{char}K\nmid n$, apply
Lemma \ref{l5.7} to get the conclusion that
$\zeta_n+\zeta^{-1}_n\in K$. Thus the condition (ii) is necessary.

For the condition (i), suppose that $n$ is even. Then
$(\bm{Z}/2\bm{Z})^2\subset D_n$. Since
$\fn{ed}_K((\bm{Z}/2\bm{Z})^2)=2$ by Theorem \ref{t4.7}, we find
that $\fn{ed}_K(D_n)\ge 2$. This is a contradiction.

Now we consider the reverse direction.

Assume that $\fn{char}K\nmid 2n$, $n$ is odd and
$\zeta_n+\zeta^{-1}_n\in K$. We find a faithful representation of
$D_n$ into $GL_2(K)$ by \cite[Lemma 7.3]{BF}. Thus we have a
faithful $D_n$-field $K(x,y)$ provided by this representation.
Define $t=\frac{y}{x}$. Then $K(t)$ is a faithful $D_n$-subfield
by \cite[Lemma 7.3]{BF} again. Thus $\fn{ed}_K(D_n)=1$.

It remain to consider the case $\fn{char}K=p\ne 2$ and $n=p$. Let
$D_p=\langle \sigma, \tau: \sigma^p=\tau^2=1,\
\tau\sigma\tau^{-1}=1\rangle$. Consider the faithful
representation $\rho:D_p\to GL_2(K)$ defined by
\[
\rho(\sigma)=\left(\begin{array}{cc}  1 & 1 \\ 0 & 1
\end{array}\right),\ \ \rho(\tau)=\left(\begin{array}{cc}  1 & 0
\\ 0 & -1  \end{array}\right).
\]

It is not difficult to show that $\fn{ed}_K(D_p)=1$ by arguments
similar to the above situation.

(2) Assume that $\fn{char}K=2$.

Suppose that $\fn{ed}_K(D_n)=1$. Then we can embed $D_n$ into
$PGL_2(K)$ by Lemma \ref{l5.3}. If $n$ is odd, then
$\zeta_n+\zeta^{-1}_n\in K$ by Lemma \ref{l5.7}. If $n$ is even,
then $n=2$ because of Lemma \ref{l5.2}. In this case, $D_n$ is
isomorphic to $(\bm{Z}/2\bm{Z})^2$. We will show that, if
$(\bm{Z}/2\bm{Z})^2$ can be embedded into $PGL_2(K)$, then $K\ne
\bm{F}_2$.

Assume to the contrary, i.e. $(\bm{Z}/2\bm{Z})^2\subset
PGL_2(\bm{F}_2)$. Note that $\fn{GL_2}(\bm{F}_2)=PGL_2(\bm{F}_2)$
and $|GL_2(\bm{F}_2)|=6$. Thus $(\bm{Z}/2\bm{Z})^2$ cannot be a
subgroup of $GL_2(\bm{F}_2)$.

It remains to prove the reverse direction. Assume that
$\fn{char}K=2$.

If $n$ is odd and $\zeta_n+\zeta^{-1}_n\in K$, then
$\fn{ed}_K(D_n)=1$ by the same arguments as in Part (1) (note that
\cite[Lemma 7.3]{BF} is valid even if $\fn{char}K=2$; the crucial
point of \cite[Lemma 7.3]{BF} is $\fn{char}K\nmid n$).

If $n=2$ and $|K|\ge 4$, choose $\alpha\in K\backslash \{0,1\}$.
Let $D_2=\langle \sigma,\tau: \sigma^2=\tau^2=1,\
\sigma\tau=\tau\sigma\rangle$. Define a faithful representation
$\rho:D_2\to GL_2(K)$ defined by
\[
\rho(\sigma)=\left(\begin{array}{cc}  1 & 1 \\ 0 & 1
\end{array}\right),\ \ \rho(\tau)=\left(\begin{array}{cc}  1 &
\alpha \\ 0 & 1  \end{array}\right).
\]

It is easy to show that $\fn{ed}_K(D_2)=1$.
\end{proof}

\begin{prop}\label{p5.9}
\r{(1)} $\fn{ed}_{\bm{F}_2}((\bm{Z}/2\bm{Z})^2)=2$.

\r{(2)} If $\fn{char}K=2$ and $|K|\ge 4$, then
$\fn{ed}_K((\bm{Z}/2\bm{Z})^2)=1$.
\end{prop}

\begin{proof}
(2) follows from Theorem \ref{t5.8}.

It remains to prove (1). From the proof of Theorem \ref{t5.8}, we
find that $\fn{ed}_{\bm{F}_2}((\bm{Z}/2\bm{Z})^2)\ge 2$. By
 Lemma \ref{l2.9} we get
$\fn{ed}_{\bm{F}_2}((\bm{Z}/2\bm{Z})^2)\leq
2\fn{ed}_{\bm{F}_2}(\bm{Z}/2\bm{Z})=2$.
\end{proof}

\begin{prop}\label{p5.10}
Let $p$ be a prime number and $K$ be any field with
$\fn{char}K=p>0$. For any positive integer $r$,
$\fn{ed}_K((\bm{Z}/p\bm{Z})^r)=1$ if and only if $[K:\bm{F}_p]\ge
r$ where $[K:\bm{F}_p]$ denotes the vector space dimension of $K$
over $\bm{F}_p$.
\end{prop}

\begin{proof}
Suppose that $[K:\bm{F}_p]\ge r$.

Let $(\bm{Z}/p\bm{Z})^r=\langle \sigma_i:\sigma^p_i=1,\
\sigma_i\sigma_j=\sigma_j\sigma_i$ for $1\le i\le r\rangle$.
Choose $\alpha_1,\alpha_2,\ldots,\alpha_r\in K$ so that
$\alpha_1,\ldots,\alpha_r$ are linearly independent over
$\bm{F}_p$. Consider the faithful representation
$\rho:(\bm{Z}/p\bm{Z})^r\to GL_2(K)$ defined by
\[
\rho(\sigma_i)=\left(\begin{array}{cc}  1 & \alpha_i \\ 0 & 1
\end{array}\right)
\]
for $1\le i \le r$. It is easy to show that
$\fn{ed}_K((\bm{Z}/p\bm{Z})^r)=1$ by similar arguments as in the
proof of Theorem \ref{t5.8}.

For the other direction, suppose that
$\fn{ed}_K((\bm{Z}/p\bm{Z})^r)=1$.

If $[K:\bm{F}_p]=\infty$, there is nothing to prove. So consider
the case that $K$ is a finite field $\bm{F}_q$ where $q=p^n$ for
some integer $n$.

Since $\fn{ed}_K((\bm{Z}/p\bm{Z})^r)=1$, we may embed
$(\bm{Z}/p\bm{Z})^r$ into $PGL_2(K)=PGL_2(\bm{F}_q)$. Let
$f:GL_2(\bm{F}_q)\to PGL_2(\bm{F}_q)$ be the canonical projection.
The group $f^{-1}((\bm{Z}/p\bm{Z})^r)$ is an extension of
$(\bm{Z}/p\bm{Z})^r$ by $\bm{F}^{\times}_q$. Since $\gcd
\{|\bm{F}^{\times}_q|,|(\bm{Z}/p\bm{Z})^r|\}=1$, we find that the
group extension splits by Schur-Zassenhaus's Theorem
\cite[p.235]{Su}. Hence $(\bm{Z}/p\bm{Z})^r$ can be embedded into
$GL_2(\bm{F}_q)$.

Since $|GL_2(\bm{F}_q)|=q(q^2-1)(q-1)$, it follows that $q$ is
divisible by $p^r$.
\end{proof}

\newpage

\end{document}